\documentclass{llncs}

\usepackage{amsmath}
\usepackage{amssymb}
\usepackage{CJK}
\usepackage{url}
\usepackage{float}
\usepackage{graphicx}
\usepackage{diagbox}
\usepackage{multirow}
\usepackage{wrapfig}
\usepackage[]{graphbox}
\newcommand{\dee}{\textrm{d}}

\begin{document}
\begin{CJK*}{UTF8}{gbsn}
\title{Evaluation of Adjoint Methods in Photoacoustic Tomography with Under-Sampled Sensors}
\titlerunning{Adjoint Photoacoustic Tomography}  
%
\author{Hongxiang Lin\inst{1} \and Takashi Azuma\inst{2} \and Mehmet Burcin Unlu\inst{3,4,5} \and Shu Takagi\inst{1}}
%
%
%
\institute{Department of Mechanical Engineering, The University of Tokyo, Tokyo, Japan\\ \email{hongxianglin@fel.t.u-tokyo.ac.jp},\\ \and Center for Disease Biology and Integrative Medicine,\\ The University of Tokyo, Tokyo, Japan \and Department of Physics, Bogazici University, Istanbul, Turkey \and Global Station for Quantum Medical Science and Engineering, Global Institution for Collaborative Research and Education, Hokkaido University, Sapporo, Japan \and Department of Radiation Oncology, Stanford University School of Medicine, Stanford, CA, USA}

\maketitle              

\begin{abstract}
Photo-Acoustic Tomography (PAT) can reconstruct a distribution of optical absorbers acting as instantaneous sound sources in subcutaneous microvasculature of a human breast. Adjoint methods for PAT, typically Time-Reversal (TR) and Back-Projection (BP), are ways to refocus time-reversed acoustic signals on sources by wave propagation from the position of sensors. TR and BP have different treatments for received signals, but they are equivalent under continuously sampling on a closed circular sensor array in two dimensions. Here, we analyze image quality with discrete under-sampled sensors in the sense of the Shannon sampling theorem. We investigate resolution and contrast of TR and BP, respectively in one source-sensor pair configuration and the frequency domain. With Hankel's asymptotic expansion to the integrands of imaging functions, our main contribution is to demonstrate that TR and BP have better performance on contrast and resolution, respectively. We also show that the integrand of TR includes additional side lobes which degrade axial resolution whereas that of BP conversely has relatively small amplitudes. Moreover, omnidirectional resolution is improved if more sensors are employed to collect the received signals. Nevertheless, for the under-sampled sensors, we propose the Truncated Back-Projection (TBP) method to enhance the contrast of BP using removing higher frequency components in the received signals. We conduct numerical experiments on the two-dimensional projected phantom model extracted from OA-Breast Database. The experiments verify our theories and show that the proposed TBP possesses better omnidirectional resolution as well as contrast compared with TR and BP with under-sampled sensors.
\end{abstract}

\section{Introduction}
Photo-Acoustic Tomography (PAT) is a prospective imaging modality that detects optical absorbers in human tissue for noninvasive diagnoses of diseases. When light is absorbed by the tissue and converted to heat, an acoustic wave is generated due to the thermoelastic expansion of the heated volume. Till now, PA microscopy, PA mammography, and PA computed tomography overcome difficulties of achieving rich optical contrast, high spatial resolution of ultrasound, as well as deep penetration depth. Nevertheless, artifacts cause image quality deterioration that significantly impacts the clinical diagnosis based on PA images. Artifacts are always concerned in a research branch called the incomplete PAT problem. In literature, PAT with a limited aperture or an inadequate broadband sensitivity at high frequency has been addressed in a sense of continuous regime. Regarding the discrete spatial sampling, deep learning based PAT can obtain high-quality images with using a training dataset~\cite{Antholzer2017}.

\begin{figure}[ht]
\vspace{-.6cm}
\begin{minipage}{.45\textwidth}
	\centering
   \includegraphics[width=0.7\textwidth,clip]{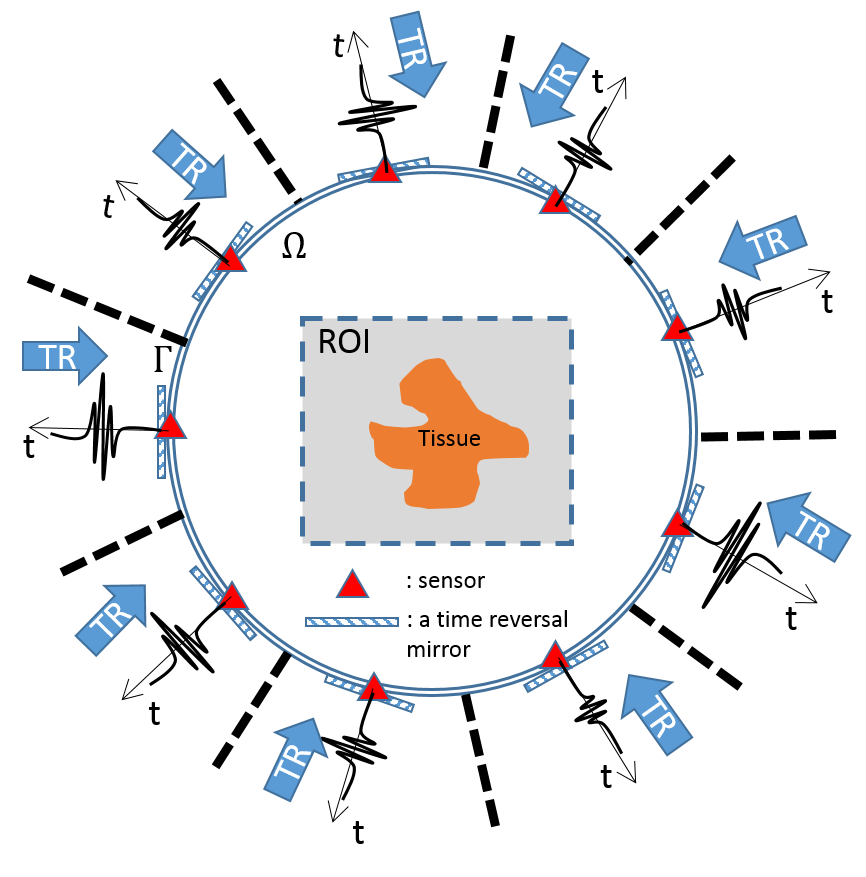}
\end{minipage}
\hspace{0.1cm}
\begin{minipage}{.45\textwidth}
	\centering
   \includegraphics[width=.5\textwidth,clip]{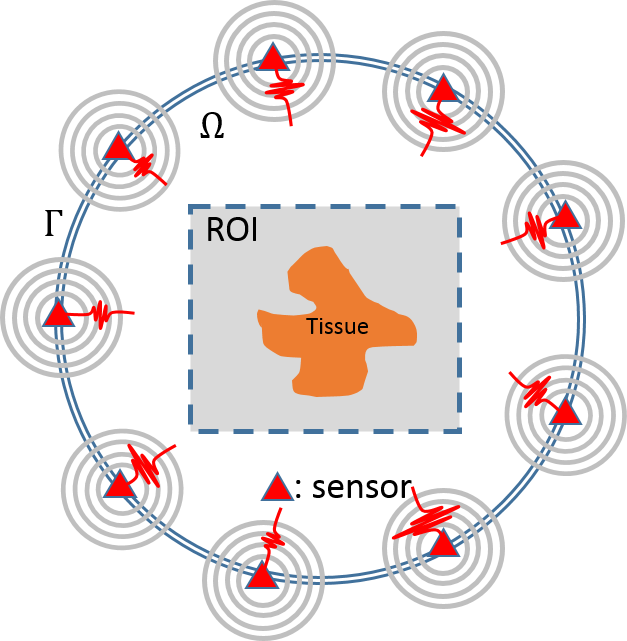}
\end{minipage}
   \label{fig:1}
   \vspace{-.5cm}
   \caption{Schematic diagrams of Time-Reversal (left) and Back-Projection (right).}
   \vspace{-.6cm}
\end{figure}

Adjoint method for photoacoustic wave propagation is a category of mathematical techniques which reverses received signals and refocuses them on source locations. In this work, we consider two typical adjoint methods -- Time-Reversal (TR) and Back-Projection (BP). As illustrated in Fig.~\ref{fig:1}, the TR method is conducted in a cavity and the reversed received waveform signals serve as a dynamic Dirichlet boundary condition. On the other hand, 
the BP method treats sensor elements as the reversing sources that retransmit circular waves modulated by reversed signals. Compared with some other iteration-based PAT methods, TR and BP possess the explicit imaging functions that illustrate the relative intensity of acoustic source distribution. The article~\cite{Ammari2011} proposes that the two methods are mathematically coincident in a continuous regime with the far-field assumption. However, the numerical study in~\cite{Arridge2016b} shows that the limited number of spatially sampled sensor elements inside a finite spatial domain may enlarge the difference between the point spread functions of TR and BP. To investigate the impact of image qualities, we quantitatively analyze resolution and contrast in the under-sampled regime, i.e., the situation where the Shannon sampling theorem is invalid~\cite{Haltmeier2016}.

Here, we focus on establishing a novel methodology to quantitatively analyze resolution and contrast of imaging functions. We unify the forms of TR and BP imaging functions in the frequency domain on one source-sensor pair configuration. Then we decompose the imaging functions as a combination of the Bessel functions since the free-space Green's function in a homogeneous medium is available. This allows us to see the contribution of main lobe for resolution and the intensity of side lobe corresponding to contrast. We also propose a straightforward variation called Truncated Back-Projection (TBP) to reconstruct the BP imaging function by removing the high-frequency components from the dataset.
%
%
\section{Mathematical Formulation}
Consider that photoacoustic wave excites from optical absorbers and propagates in a two-dimensional lossless homogeneous medium with a speed of sound of $c_0$; see Fig.~\ref{fig:1}. The circular boundary $\Gamma$ of a radius of $R$ is composed by $N$ point-like sensor elements located at $\mathbf{y}_n\ (n=1,2,\cdots,N)$ with equispaced arguments. The photoacoustic signals denoted by $g(\mathbf{y}_n, t)$ are semi-discretized in a time interval $[0, T]$ and sensor elements $\{\mathbf{y}_n\}$ on $\Gamma$. The purpose of PAT is to recover the initial pressure distribution $p_0$ in a region of interest (ROI) given the photoacoustic measurements $g(\mathbf{y}_n, t)$ on $\Gamma\times [0,T]$.

Adjoint methods are derived to approximate $p_0$ through different treatments for the reversed signals $g(\mathbf{y}_n, T-t)$ as well as refocusing the reversed wavefield on $p_0$ at the terminal time $T$.
Together with the Green's representation theorem, the boundary-condition treatment for $g(\mathbf{y}_n, T-t)$ yields the semi-discrete TR imaging function:
\begin{equation}\label{eq:TR,semi}
I_{TR}^{<N>}(\mathbf{x}) = h_N\sum_{n=1}^N\int_0^T \frac{\partial G_d}{\partial \mathbf{\nu_y}}(\mathbf{x},T|\mathbf{y}=\mathbf{y}_n,t) g(\mathbf{y}_n,T-t)dt+O(h_N^2),
\end{equation}
where $G_d$ is the Dirichlet Green's function of wave equation, $\mathbf{\nu_y}$ an outward unit normal vector at $\mathbf{y}$, and $h_N = 2\pi R/N$ a step size along $\Gamma$. 
On the other hand, the reversing-source treatment for $g(\mathbf{y}_n, T-t)$ infers the semi-discrete BP imaging function:
\begin{equation}\label{eq:MTR,semi}
I_{BP}^{<N>}(\mathbf{x}) = \frac{h_N}{c_0}\sum_{n=1}^N\int_0^T\frac{\partial G_0}{\partial t}(\mathbf{x},\tau| \mathbf{y}_n,t=T)g(\mathbf{y}_n,T-\tau)d\tau + O(h_N^2),
\end{equation}
where $G_0$ is the free-space Green's function of wave equation. The derivation of Eqs.~\ref{eq:TR,semi} and~\ref{eq:MTR,semi} is referred to~\cite{Arridge2016b} or the supplementary material.

The mathematical analysis is henceforth conducted in the frequency domain. The Fourier transform of a time-history function $f(t)$ is defined as $\hat{f}(\omega) = \int_{-\infty}^{+\infty}f(t)e^{i\omega t}dt$ where $\omega$ is an angular frequency and the hat $\wedge$ denotes the Fourier transform. Using the Parseval's identity (see the supplementary material) to Eqs.~\ref{eq:TR,semi} and~\ref{eq:MTR,semi}, we write out the frequency-domain expressions for TR and BP:
\begin{eqnarray}
I_{TR}^{<N>}(\mathbf{x}) = \frac{h_N}{2\pi}\mathbf{Re}\sum_{n=1}^N\left\{\int_{-\infty}^{+\infty}\widehat{\frac{\partial G_d}{\partial\mathbf{\nu_y}}}(\mathbf{x},\mathbf{y}=\mathbf{y}_n,\omega)\overline{\widehat{g}(\mathbf{y}_n,\omega)}d\omega\right\} + O(h_N^2), \label{eq:TR,parseval} \\
I_{BP}^{<N>}(\mathbf{x}) = -\frac{h_N}{2\pi c_0}\mathbf{Re}\sum_{n=1}^N\left\{\int_{-\infty}^{+\infty} i\omega\widehat{G_0}(\mathbf{x},\mathbf{y}_n,\omega)\overline{\widehat{g}(\mathbf{y}_n,\omega)}d\omega\right\} + O(h_N^2), \label{eq:MTR,parseval}
\end{eqnarray}
where $\mathbf{Re}$ denotes the real part of a complex value and the overline denotes complex conjugate. The free-space Green's function in the frequency domain is written as $\widehat{G_0}(\mathbf{x},\mathbf{y},\omega) = \frac{i}{4}H_0^{(1)}\left(\frac{\omega}{c_0}|\mathbf{x}-\mathbf{y}|\right)$
where $H_0^{(1)}$ is a zeroth-order Hankel function of the first kind.

\section{Image Quality Analysis}
One source-sensor pair configuration is considered to characterize the image quality of TR and BP images. 
There is only one acoustic source located at $\mathbf{a}$ in the cavity $\Omega$. The sensor at $\mathbf{y}$ receiving the single-source waveform signal satisfies the frequency-domain expression: 
\begin{equation}\label{eq:pt_response}
\hat g(\mathbf{y},\omega)=-i\omega F(\omega)\widehat{G_0}(\mathbf{y},\mathbf{a},\omega),
\end{equation}
where $F(\omega)$ is a real function of the $\omega_{\max}$-bandlimited spectrum.

\subsection{Expansion of Imaging Functions}
Substitute Eq.~\ref{eq:pt_response} into Eqs.~\ref{eq:TR,parseval} and~\ref{eq:MTR,parseval} first, which yields a unified imaging function for the source-sensor pair $(\mathbf{a}, \mathbf{y})$ configuration:
\begin{equation}\label{eq:img_fun,integrand}
I_j^{<1>}(\mathbf{x}) = \frac{h_1}{2\pi}\int_{-\infty}^{+\infty}F(\omega)\ \mathrm{Re} [K_j(\mathbf{x},\omega)]d\omega + O(h_1^2),\ j=\textrm{TR, BP}
\end{equation}
where the integrands are specified as $
K_{TR}(\mathbf{x},\omega) = i\omega\widehat{\frac{\partial G_d}{\partial\mathbf{\nu_y}}}(\mathbf{x},\mathbf{y},\omega)\overline{\widehat{G_0}(\mathbf{y},\mathbf{a},\omega)}, \label{eq:TR,integrand}$
and
$K_{BP}(\mathbf{x},\omega) = \frac{\omega^2}{c_0}\widehat{G_0}(\mathbf{x},\mathbf{y},\omega)\overline{\widehat{G_0}(\mathbf{y},\mathbf{a},\omega)}. \label{eq:MTR,integrand}$

\setcounter{footnote}{0}
By canceling factors, we reduce $K_{TR}$ and $K_{BP}$ to $\widetilde{K_{TR}}$ and $\widetilde{K_{BP}}$, respectively,
such that they share an identical main lobe.\footnote{In specific, $\widetilde{K_{TR}} = -32\pi c_0^3\frac{|\mathbf{y}-\mathbf{a}|}{|\mathbf{y}|}K_{TR}$ and $\widetilde{K_{BP}} = 32\pi c_0^3K_{BP}$.} Employing Hankel's asymptotic expansion~\cite{Abramowitz1966}, we write out the BP imaging function $I_{BP}^{<1>}$ and the discrepancy function $\Delta I^{<1>}$ between TR and BP in proportion to the integration of $\widetilde{K_{BP}}$ and $\widetilde{K_{TR}}-\widetilde{K_{BP}}$ over the angular frequency domain respectively:
\begin{align}
I_{BP}^{<1>}(\mathbf{x}) &\propto \int_{-\infty}^{+\infty}F(\omega)\ \mathrm{Re}[\widetilde{K_{BP}}](\mathbf{x},\omega)\dee\omega\notag \\ 
&\approx \int_{-\infty}^{+\infty}F(\omega)\left[\frac{2\omega c_0}{\pi|\mathbf{y}-\mathbf{a}|}J_0(d_{\mathbf{x},\mathbf{a}}^{(\omega)})+ O((d_{\mathbf{x},\mathbf{a}}^{(\omega)})^3)\right]\dee\omega, \label{eq:mtr,integrand,expansion,psf} \\
\Delta I^{<1>}(\mathbf{x}) &\propto \int_{-\infty}^{+\infty}F(\omega)\ \mathrm{Re}[\widetilde{K_{TR}}-\widetilde{K_{BP}}](\mathbf{x},\omega)\dee\omega\notag \\
&\approx e^{i\Theta}\int_{-\infty}^{+\infty}F(\omega)\left[\frac{2c_0^2}{\pi|\mathbf{y}-\mathbf{a}|^2}J_1(d_{\mathbf{x},\mathbf{a}}^{(\omega)})+ O((d_{\mathbf{x},\mathbf{a}}^{(\omega)})^4)\right]\dee\omega  \label{eq:diff,integrand,expansion,psf}
\end{align}
where $d_{\mathbf{x},\mathbf{a}}^{(\omega)} = \frac{\omega}{c_0}|\mathbf{x}-\mathbf{a}|$. $\Theta$ is the angle corresponding to the opposite side $|\mathbf{x}-\mathbf{y}|$ of the triangle formed by the points $\mathbf{a}$, $\mathbf{x}$ and $\mathbf{y}$. See the derivation in the supplementary material. If $\omega_c$ is the angular center frequency and $\Theta = 0$ or $\pi$, it yields that the axial pattern of BP is approximately the zeroth-order Bessel function $J_0(d_{\mathbf{x},\mathbf{a}}^{(\omega_c)})$ while that of TR has the same main lobe as BP plus a side lobe of the first-order Bessel function $J_1(d_{\mathbf{x},\mathbf{a}}^{(\omega_c)})$.

\subsection{Resolution Analysis}\label{sec:3.3}
Axial resolution is quantified by Full Width at Half Maximum (FWHM). Based on the axial pattern in Eq.~\ref{eq:mtr,integrand,expansion,psf}, since $J_0(\xi)$ has a maximum at $\xi = 0$ and a half maximum at approximately $\xi = 1.5$, we have FWHM of BP: $W_{\mathrm{BP}}^{\mathrm{FWHM}} \approx 2\times 1.5/(2\pi/\lambda_c) \approx 0.48\lambda_c$ where $\lambda_c$ is the wavelength corresponding to the center frequency. Similarly, since the half maximum values of $J_1(\xi)$ are located at $\xi = 0.6$ and $3.1$, we have FWHM of the side lobe shown in Eq.~\ref{eq:diff,integrand,expansion,psf}: $W_{\Delta I}^{\mathrm{FWHM}} \approx (3.1-0.6)/(2\pi/\lambda_c)\approx 0.40\lambda_c$. Additionally, for both TR and BP, the lateral resolution degrades by noting that the radial transmission of wavelet implied in Eq.~\ref{eq:img_fun,integrand} leads to artifact of an arc pattern.

\begin{figure}[ht]
\vspace{-.6cm}
\begin{minipage}[t]{.3\textwidth}
	\centering
   \includegraphics[height=0.15\textheight,clip]{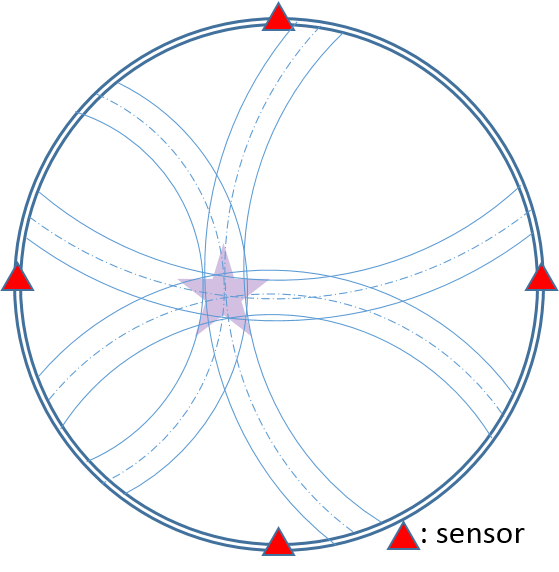} 
   \vspace{-.3cm}  
   \caption{The sketch of the adjoint wavefield synthesized by four equispaced-sensor signals refocused on the source (star).}
   \label{fig:2}
\end{minipage}
\hspace{0.2cm}
\begin{minipage}[t]{.68\textwidth}
	\centering
   \includegraphics[height=0.15\textheight,clip]{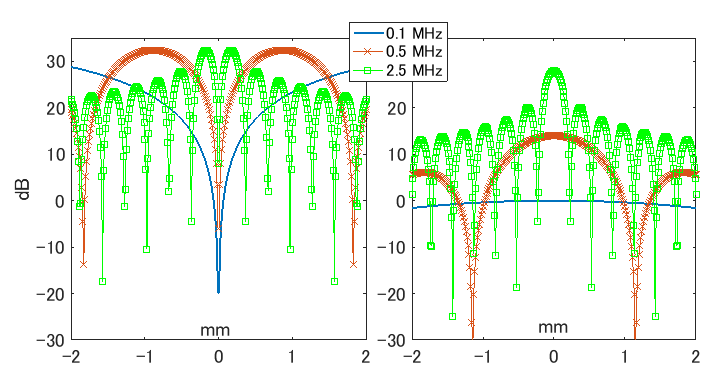}
   \vspace{-.3cm}
   \caption{The intensity level of $\widetilde{K_{TR}}-\widetilde{K_{BP}}$ (left) and $\widetilde{K_{BP}}$ (right) on one source-sensor pair configuration along the axial direction. The source $\mathbf{a}$ is placed at the origin. The reference value is given by $\widetilde{K_{BP}}$ at the central point and at a frequency of $0.1$ MHz.}
     \label{fig:3}
\end{minipage}
\vspace{-.6cm}
\end{figure}

For adjoint methods, the omni-directional resolution of the source point $\mathbf{a}$ can be extended from the axial one through superposition of adjoint wavelets. As shown in Fig.~\ref{fig:2}, the four wavelets serve as carriers of the same FWHM information oriented from the different directions. In morphology, they partially overlap and form a polygon-like spot approximating a circle. Moreover, the superposed wavefield mitigates the artifact since the amplification at $\mathbf{a}$ significantly inhibits the intensity level of synthesized wavefield in other pixels where there are no sources. 


\subsection{Contrast Analysis and Truncated Back-Projection}
The contrast is discussed in the regime of an under-sampled spatial grid; the full-sampling condition is referred to~\cite{Haltmeier2016}. Figure~\ref{fig:3} shows the profiles of $\widetilde{K_{TR}}-\widetilde{K_{BP}}$ and $\widetilde{K_{BP}}$ along the axial direction with respect to three typical frequencies. The side lobe of the imaging function of TR dominates due to the larger intensity in a range of frequency.
Therefore, for TR, the intensity level in a prescribed neighborhood of  the source spot does not fiercely alter with the change of the spatial sampling rate. Conversely, regrading BP, since we often interpolate the circular-band wavefield on a coarse Cartesian grid, the maximum may occasionally be selected out but the rest of the interpolated grid points have no significantly large values. Their contrast deteriorates after normalization is conducted.

To solve this, we propose the Truncated Back-Projection (TBP) method by means of only exploiting the low frequency components in Eq.~\ref{eq:MTR,parseval} given by
\begin{equation}\label{eq:TBP_parseval}
I_{TBP}^{<N,\mu>}(\mathbf{x}) = -\frac{h_N}{2\pi c_0}\mathbf{Re}\sum_{n=1}^N\left\{\int_{|\omega|<\mu} i\omega\widehat{G_0}(\mathbf{x},\mathbf{y},\omega)\overline{\widehat{g}(\mathbf{y},\omega)}d\omega\right\} + O(h_N^2),
\end{equation}
where $\mu$ is the truncated bound of angular frequency. If we select a $\mu$ much smaller than the upper bound $\omega_{\max}$,
more large values adjacent to the maximum of $J_0(d_{\mathbf{x},\mathbf{a}}^{(\mu)})$ are attainable in the coarse grid
since the larger FWHM proportional to the wavelength has capacity of containing more grid points. We recommend to set $\mu=2\pi M/(T_M\cdot PPW)$ where $M$ is the number of gird points per side of ROI, $T_M$ acquisition time in ROI, and $PPW$ the number of grid points per wavelength indicating the coarseness of grid.

\section{Numerical Experiments}
We carry out two numerical experiments with a single source and a breast vasculature phantom model used as initial pressure distributions. The acoustic measurement datasets are synthesized by the K-wave toolbox, a photoacoustic Matlab simulator using the k-space pseudo spectral method~\cite{Treeby2010}. To avoid the inverse crime, the measurement datasets are generated on a fine Cartesian grid but the adjoint methods are conducted on a coarse one. Their source codes are available at github for reproduction of results.\footnote{\url{https://github.com/hongxiangharry/AdjointPAT}}

\subsection{Single Source-Sensor Pair Reconstruction}\label{sec:4.1}
\begin{figure}[ht]
   \vspace{-.6cm}
   \centering
   \includegraphics[width=1\textwidth,clip]{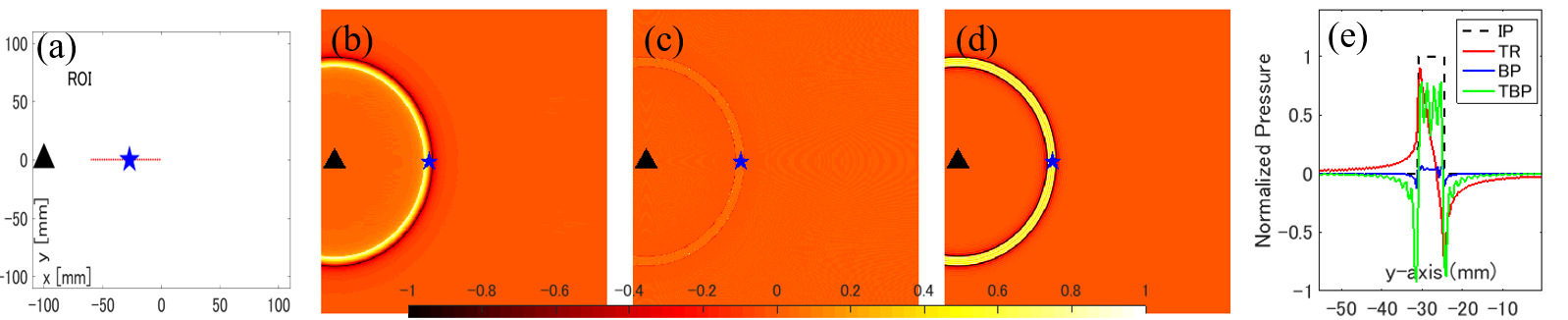}
   \vspace{-.5cm}
   \caption{(a) One source-sensor pair configuration with the source at $(-12.5, 0)$ [mm] and the sensor at $(-100, 0)$ [mm]. The adjoint wavefields reconstructed by (b) TR, (c) BP, and (d) TBP. (e) The normalized pressure distributions of the initial pressure (IP) distribution, TR, BP, and TBP profiles along the red dot line in (a). The black triangles and blue stars represent sensors and sources, respectively.}
   \label{fig:5}
   \vspace{-1cm}
\end{figure}

Figure~\ref{fig:5} shows the circular-arc wavefields reconstructed by the adjoint methods. All the wavefields are normalized by the maximum absolute values. In terms of axial resolution, the values of FWHM for TR, BP, and TBP are $2.67$, $6.19$, and $5.58$ mm shown in Fig.~\ref{fig:5}(e). The maximum values, used to assess the contrast, for TR, BP, and TBP are $0.90$, $0.089$, and $0.78$. We observe that the axial profile of the TR image has a fat-tail distribution because of relatively strong side lobe.

In terms of spatial sampling, the maximum is sparsely selected out from a coarse grid. It may lead to the deterioration of contrast in the BP image. The numerical experiment is conducted on a grid of $512$-by-$512$ points, coarse for $PPW = 0.89$, provided $\omega_{\max} = 2.73\times 10^7\ \mathrm{rad}\cdot\mathrm{s}^{-1}$ and $T_M=1.33\times 10^{-4}$ s. It requires at least a $2500$-by-$2500$ grid to achieve the regime of full sampling, i.e. $PPW = 4.32$. The truncated bound $\mu$ limited to $6.82\times 10^6\ \mathrm{rad}\cdot\mathrm{s}^{-1}$ for TBP fulfills the need that the $512$-by-$512$ grid satisfies $PPW = 4.32$ although its point spread function is oscillating.

\subsection{Breast Phantom Reconstruction}
\begin{figure}[ht]
\vspace{-.6cm}
\centering
\begin{tabular}{c|l}
Ele. & $\ $Configuration$\qquad\qquad$TR$\qquad\qquad\quad\ $BP$\qquad\qquad\quad\ $TBP$\qquad\qquad\quad\ $Profile \\
\hline
$16$ & \includegraphics[align=c,width=.95\textwidth,clip]{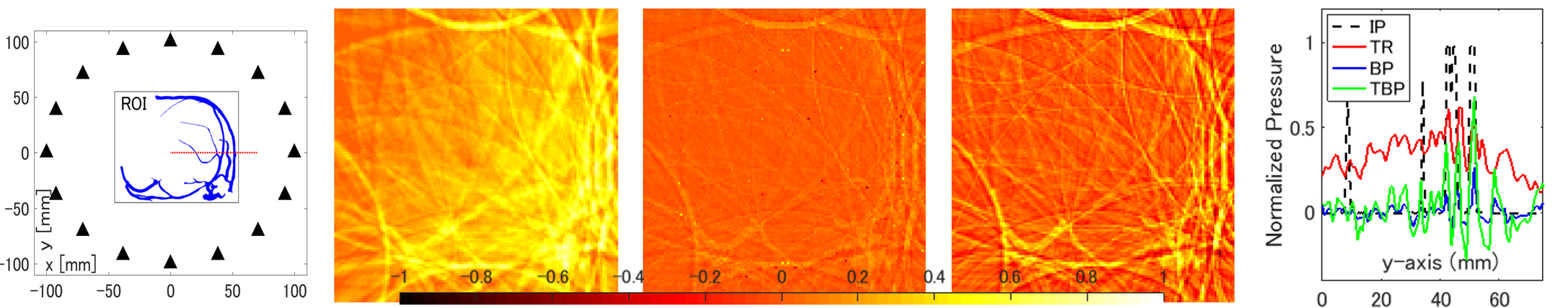}\\ 
$64$ & \includegraphics[align=c,width=.95\textwidth,clip]{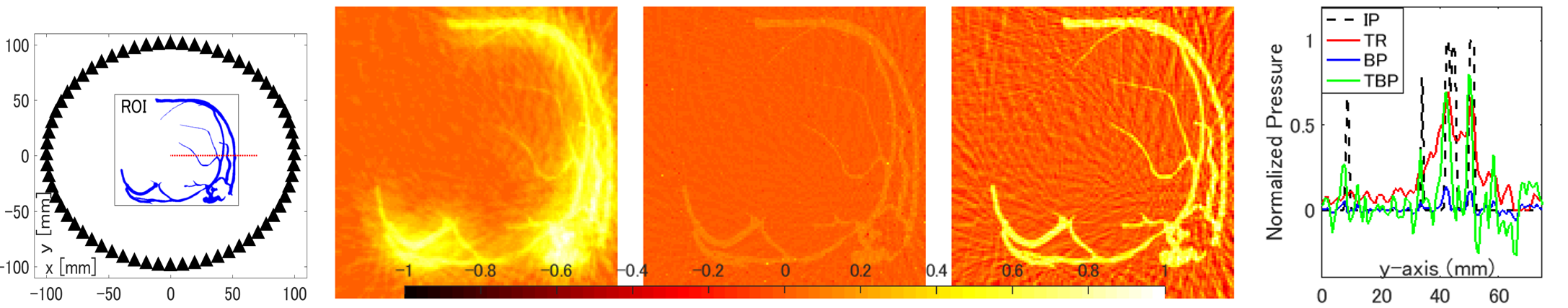}\\
$256$ & \includegraphics[align=c,width=.95\textwidth,clip]{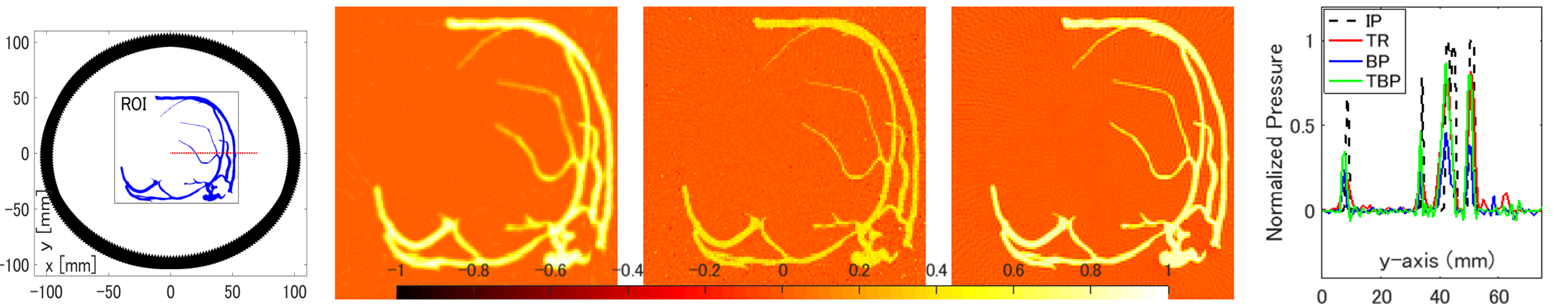}\\
\end{tabular}
\vspace{-.3cm}
     \caption{The breast vasculature images reconstructed by the adjoint methods with the different sampling of sensor elements. ROI is a square centered at the origin with an area of $100\times 100\ \textrm{mm}^2$. The reconstructions are carried out on three arrays with equispaced sensor elements (ele.). In the configuration, the black triangles and blue lines represent sensors and sources, respectively. The images reconstructed by TR, BP, and TBP are shown in the middle three columns. The last column shows the profiles along the red dot line in the first column.} 
       \label{fig:6} 
       \vspace{-.6cm}
\end{figure}

The breast phantom is extracted from OA-Breast Database~\cite{Lou2017} and projected to two dimensions. The reconstructed vasculature images with respect to $16$, $64$, and $256$ equispaced sensor elements are shown in Fig.~\ref{fig:6}, which corresponds to the under sampling, the critical-condition sampling, and the full sampling based on Shannon sampling theorem. The more sensor elements we utilize, the better resolution and contrast will achieve for TR, BP, and TBP, which validates the summary in Sect.~\ref{sec:3.3}. We demonstrate that whenever the sampling criterion is selected, TR and BP have advantages of contrast and axial resolution, respectively. Moreover, a $256$-by-$256$ grid of $PPW=1.32$ are used to validate TBP, provided $\omega_{\max} = 1.82\times 10^7\ \mathrm{rad}\cdot\mathrm{s}^{-1}$ and $T_M=6.67\times 10^{-5}$ s. The parameter $\mu = 6.82\times 10^6\ \mathrm{rad}\cdot\mathrm{s}^{-1}$ assures high contrast of the TBP image even in the situation of the under-sampled sensors.

\section{Conclusions}
In this work, we demonstrate that TR and BP possess high contrast and high axial resolution with under-sampled sensors, respectively. Asymptotic expansion technique helps to mathematically specify the intrinsic behaviors of TR and BP in the frequency domain. We propose the TBP method to compensate the contrast issue in the situation of spatial coarse grid. Although all of the analyses and numerical tests are presented in two dimensions, the methodology is possible to be correspondingly extended to high dimensions with the irregular geometry of sensor array and the complex structure of medium~\cite{Borcea2003}.



\end{CJK*}
\end{document}